\documentclass[10pt]{amsart}
\usepackage{graphicx}
\usepackage[USenglish]{babel}
\usepackage{latexsym}
\usepackage{amsmath}
\usepackage{amsfonts}
\usepackage{amssymb}
\usepackage{esint}
\renewcommand{\a }{\alpha }

\renewcommand{\d}{\delta }
\newcommand{\D }{\Delta }

\renewcommand{\l }{\lambda }
\renewcommand{\L }{\Lambda }

\newcommand{\n }{\nabla }
\newcommand{\var }{\varphi }
\newcommand{\rh }{\rho }

\newcommand{\s }{\sigma }
\newcommand{\Sg }{\Sigma}

\renewcommand{\th }{\theta }

\renewcommand{\O }{\Omega }

\newcommand{\pa }{\partial}

\newcommand{\ov}{\overline}

\newcommand{\wtilde }{\widetilde}
\newcommand{\wk}{\rightharpoonup}
\newcommand{\be}{\begin{equation}}
\newcommand{\ee}{\end{equation}}
\newenvironment{pf}{\noindent{\sc Proof}.\enspace}{\rule{2mm}{2mm}\medskip}
\newenvironment{pfn}{\noindent{\sc Proof}}{\rule{2mm}{2mm}\medskip}

\newcommand{\R}{\mathbb{R}}

\newcommand{\N}{\mathbb{N}}

\newcommand{\dis}{\displaystyle}

\newtheorem{lem}{Lemma}[section]
\newtheorem{pro}[lem]{Proposition}
\newtheorem{thm}[lem]{Theorem}
\newtheorem{rem}[lem]{Remark}
\newtheorem{cor}[lem]{Corollary}
\newtheorem{df}[lem]{Definition}

    \newenvironment{dedication}
        {\vspace{0ex}\begin{quotation}\begin{center}\begin{em}}
        {\par\end{em}\end{center}\end{quotation}}

\begin{document}

\title[A Moser-Trudinger inequality for the singular Toda system] {A Moser-Trudinger inequality for the singular Toda system}

\author{Luca Battaglia$^{(1)}$, Andrea Malchiodi$^{(2)}$}

\address{$^{(1)}$ SISSA, via Bonomea 265, 34136 Trieste (Italy).}

\address{$^{(2)}$ University of Warwick, Mathematics Institute, Zeeman Building, Coventry CV4 7AL and SISSA, via Bonomea 265, 34136 Trieste (Italy).}

\thanks{The authors are supported by the FIRB project {\em Analysis and Beyond} and by the PRIN {\em Variational Methods and Nonlinear PDEs}. L.B. acknowledges support from the Mathematics Department at the University of Warwick.}

\email{lbatta@sissa.it, A.Malchiodi@warwick.ac.uk, malchiod@sissa.it}

\begin{abstract}
In this paper we prove a sharp version of the Moser-Trudinger inequality for the Euler-Lagrange functional of a 
singular Toda system, motivated by the study of models in Chern-Simons theory. Our result extends those in 
\cite{chenwx} and \cite{tro} for the scalar case, as well as that in \cite{jw} for the regular Toda system. We expect this inequality 
to be a  basic tool to attack variationally the  existence problem under general assumptions. 
\end{abstract}

\maketitle

    \begin{dedication}
    Dedicated to Neil Trudinger with admiration
    \end{dedication}

\section{Introduction}

\noindent The  Moser-Trudinger inequality yields exponential-type embeddings of Sobolev functions in 
critical dimension. On a compact closed surface $\Sigma$ the space $H^1(\Sigma)$ embeds compactly into 
every $L^p(\Sigma)$ for any real $p > 1$: at a more refined level, due to the seminal works \cite{tru} 
and \cite{mos} one has the inequality 
\begin{equation}\label{eq:mt}
 16 \pi\log\int_{\Sigma}e^{u-\overline{u}}\,dV_{g}\leq\int_{\Sigma}\left|\nabla u\right|^{2}\,dV_{g}+C;
 \qquad \quad u \in H^1(\Sigma),
\end{equation}
where $C$ is a constant depending only on $\Sigma$ and its metric $g$, and where $\overline{u}$ stands for the 
average of $u$ on the surface. 

Inequality \eqref{eq:mt} has been proven to be fundamental in several contexts such as the Gaussian curvature 
prescription problem (\cite{aub}, \cite{cy0}, \cite{cysing}),  mean field equations in fluid dynamics  (\cite{djlw}, 
\cite{dja}) and models in theoretical physics
(\cite{tarl}, \cite{yyang}). To give an example, 
considering a conformal change of metric 
of the form $\wtilde{g} = e^{w} g$, the Gaussian curvature of $\Sigma$ transforms according to the law 
\begin{equation}\label{eq:gc}
- \Delta w + 2 K_g = 2 K_{\wtilde{g}} e^{w}. 
\end{equation}
If one wishes to prescribe the Gaussian curvature $K_{\wtilde{g}}$ as a given function $K(x)$, then solutions 
to the problem can be found as critical points of the functional 
$$
  I(u) := \int_{\Sigma}  |\nabla u|^2 dV_g + \int_\Sigma K_g u \, dV_g - 
  \left( \int_\Sigma K_g dV_g \right) \log \left( \int_{\Sigma} K \, e^{u} dV_g \right). 
$$
By means of \eqref{eq:mt} one can then control the last term in the functional by means of the Dirichlet energy.

More recent versions  of \eqref{eq:mt} include exponential terms with power-type weights, which are motivated 
by the study of {\em singular Liouville equations}.  For example, given points $p_1, \dots, p_m \in \Sigma$,  
weights $\a_1, \dots, \a_m > - 1$, and a smooth positive function $h(x)$,  a solution of the equation 
\begin{equation}\label{eq:gcsing}
- \Delta w + 2 K_g = 2 h \,  e^{w} - 4 \pi \sum_{j=1}^{m} \alpha_j \delta_{p_j} 
\end{equation}
yields a conformal metric $\wtilde{g} = e^w g$ with Gaussian curvature $h$ on $\Sigma 
\setminus \{ p_1, \dots, p_m\}$ and with a conical singularity at $p_j$ with opening angle $2 \pi (1 + \alpha_j)$. 

By the substitution
\begin{equation}\label{eq:changesing}
    w(x) \mapsto w(x) + 4\pi \sum_{j=1}^m \alpha_{i} G_{p_j}(x),\qquad h(x)\mapsto \wtilde{h}(x)= h(x) e^{-4\pi \sum_{j=1}^m \alpha_{j} G_{p_j}(x)}, 
\end{equation}
\eqref{eq:gc} transforms into an equation of the form
\begin{equation}\label{eq:gcsingreg}
- \Delta w + 2 \wtilde{f} = 2 \wtilde{h} \, e^{w} 
\end{equation}
where $\wtilde{f}(x)$ is a smooth function and where  
\begin{equation} \label{eq:h}
\wtilde{h} > 0 \quad \hbox{ on } \Sigma \setminus \{p_1, \dots,
p_{m}\}; \qquad \quad
  \wtilde{h}(x) \simeq d(x,p_j)^{2 \a_{j}} \mbox{ near } p_j.
\end{equation}
Although \eqref{eq:gcsing} and \eqref{eq:gcsingreg} are perfectly equivalent, the advantage 
of the latter compared to the former is that the singular structure is absorbed into the factor $\wtilde{h}$, 
which endows the problem with a variational structure.  Similarly to \eqref{eq:gc}, 
solutions to \eqref{eq:gcsingreg} can be found as critical points of the functional 
$$
  \wtilde {I}(u) := \int_{\Sigma}  |\nabla u|^2 dV_g + \int_\Sigma \wtilde{f} \, u \, dV_g - 
 \ov{K} \log \left( \int_{\Sigma} \wtilde {h} \,  e^{u} dV_g \right),
$$
where $\ov{K} = 2 \pi \chi(\Sigma) + 2 \pi \sum_{j=1}^m \alpha_j$ is a constant determined by the 
Gauss-Bonnet formula.

The singular weight $\wtilde{h}$ has indeed an effect on the optimal constant in the 
corresponding Moser-Trudinger type inequality. In \cite{chen}, \cite{tro} (see also \cite{cysing} for conical domains) 
it was shown that 
\begin{equation}\label{eq:mtsing}
 16 \pi 
 \min \left\{ 1, 1 + \min_j \alpha_j \right\}\log\int_{\Sigma} \wtilde {h} e^{u-\overline{u}}\,dV_{g}\leq\int_{\Sigma}\left|\nabla u\right|^{2}\,dV_{g}+C;
 \qquad \quad u \in H^1(\Sigma). 
\end{equation}
Notice that, if at least one of the $\alpha_j$'s is negative, say $\alpha_{\ov{j}}$, the constant gets worse, 
as $\wtilde{h}$ blows-up at $p_{\ov{j}}$. On the other hand when all the weights are positive the costant does not improve: 
this can be easily seen by the following consideration. The sharpness of the Moser-Trudinger constant $\frac{1}{16\pi}$ 
can be obtained using the test function 
\begin{equation}\label{eq:phil}
  \var_{\l,x}(y) = \log \frac{\l^2}{\left( 1 + \lambda^2  d(x,y)^2\right)^4}; \qquad 
  \quad x \in \Sigma, \l > 0, 
\end{equation}
which makes the two sides of \eqref{eq:mt} diverge at the same rate. As the conformal volume $e^{\varphi_{\lambda,x}}$ 
concentrates at $x$ as $\l \to + \infty$, there would be no effect from the vanishing of $\wtilde{h}$ if $x$ is 
a regular point. We also refer to \cite{cia}, \cite{fm} for more general optimal inequalities on singular 
measure spaces. 

Inequality \eqref{eq:mtsing} has been useful in finding constant curvature metrics when prescribing 
conical singularities as it might yield  global minima of $\wtilde{I}$, see \cite{tro}, \cite{carma}, as well as in studying 
general singular mean field equations like 
\begin{equation}\label{eq:singmf}
- \Delta w + 2 f = 2 \rho h \, e^{w} - 4 \pi \sum_{j=1}^{m} \alpha_j \delta_{p_j},  
\end{equation}
where $f, h$ are smooth functions, $h$ positive, and $\rho$ is a real parameter, see \cite{bm}, \cite{bdm}, \cite{malrui} (see also 
\cite{clin}, \cite{cl3} for a non-variational approach to \eqref{eq:singmf}). 

Singular Liouville equations have a  role in fluid dynamics, see \cite{ty}, as well as in the study of Electroweak theory or 
abelian Chern-Simons vortices, see \cite{tarl}, \cite{yyang}. For the latter cases, 
singular points represent zeroes of the scalar wave function involved in the model.

\

\noindent The goal 
of this paper is to prove a sharp inequality related to a {\em singular Toda system} arising in 
Chern-Simons theory,  which represents a non-abelian counterpart of 
\eqref{eq:singmf}. Specifically, we consider the following system 
\begin{equation}\label{eq:e-1}
  \left\{
      \begin{array}{ll}
        - \D u_1 = 2 \rho_1 \left( \frac{h_1 e^{u_1}}{\int_\Sg
      h_1 e^{u_1} dV_g} - 1 \right) - \rho_2 \left( \frac{h_2 e^{u_2}}{\int_\Sg
      h_2 e^{u_2} dV_g} - 1 \right) - 4 \pi \sum_{j=1}^{m} \a_{1,j} (\d_{p_j}-1), \\
       - \D u_2 = 2 \rho_2 \left( \frac{h_2 e^{u_2}}{\int_\Sg
      h_2 e^{u_2} dV_g} - 1 \right) - \rho_1 \left( \frac{h_1 e^{u_1}}{\int_\Sg
      h_1 e^{u_1} dV_g} - 1 \right) - 4 \pi \sum_{j=1}^{m} \a_{2,j} (\d_{p_j}-1), &
      \end{array}
  \right. \end{equation}
where $h_1, h_2$ are smooth positive functions  on $\Sigma$, and the coefficients $\alpha_{i,j}$ 
are larger than $-1$.

While abelian Chern-Simons vortices have been quite studied for some time,
see e.g. \cite{cay}, \cite{ci}, \cite{ntcv99}, \cite{sy}, \cite{tar96}, the treatment of the non-abelian
case is more recent, see e.g. \cite{dunne1995self}, \cite{kaolee}, \cite{lee},
\cite{noltar2000}.

With a change of variable similar to \eqref{eq:changesing} the latter problem transforms into 
\begin{equation} \label{eq:todaregul}
\left\{
    \begin{array}{ll}
      - \D u_1 = 2 \rho_1 \left( \frac{\wtilde{h}_1 e^{u_1}}{\int_\Sg
    \wtilde{h}_1 e^{u_1} dV_g} - 1 \right) - \rho_2 \left( \frac{\wtilde{h}_2 e^{u_2}}{\int_\Sg
    \wtilde{h}_2 e^{u_2} dV_g} - 1 \right), \\
     - \D u_2 = 2 \rho_2 \left( \frac{\wtilde{h}_2 e^{u_2}}{\int_\Sg
    \wtilde{h}_2 e^{u_2} dV_g} - 1 \right) - \rho_1 \left( \frac{\wtilde{h}_1 e^{u_1}}{\int_\Sg
    \wtilde{h}_1 e^{u_1} dV_g} - 1 \right), &
    \end{array}
  \right.
\end{equation}
where the functions $\wtilde{h}_i$ satisfy
\begin{equation} \label{eq:hi}
\wtilde{h}_i > 0 \quad \hbox{ on } \quad \Sigma \setminus \{p_1, \dots,
p_{m}\}; \qquad 
  \wtilde{h}_i(x) \simeq d(x,p_j)^{2 \a_{i,j}} \mbox{ near } p_j, \qquad \quad i=1,\ 2.
\end{equation}
As for the scalar case one gains the variational structure, with Euler-Lagrange functional  
\begin{equation} \label{funzionale}
J_{\rho}(u_1, u_2) = \int_\Sg Q(u_1,u_2)\, dV_g  + \sum_{i=1}^2
\rho_i \left ( \int_\Sg u_i dV_g - \log \int_\Sg \wtilde{h}_i e^{u_i} dV_g
\right ),
\end{equation}
 where $Q(u_1,u_2)$ is defined as:
\begin{equation}\label{eq:QQ}
Q(u_1,u_2) = \frac{1}{3} \left ( |\n u_1|^2 + |\n u_2|^2 + \n u_1 \cdot
\n u_2\right ).
\end{equation}
Concerning Liouville systems with no singularites, some sharp inequalities were proven in \cite{csw}, \cite{wang99} when the matrix of 
coefficients of the exponential terms is non-negative. For the regular Toda system instead a sharp inequality 
was found in \cite{jw}, where it was shown that 
\begin{equation}\label{eq:jw}
   4 \pi\sum_{i=1}^2  \log\int_{\Sigma}e^{u_i-\overline{u_i}}\,dV_{g}
   \leq\int_{\Sigma}Q(u_1, u_2)\,dV_{g}+C;
    \qquad \quad u \in H^1(\Sigma). 
\end{equation}
Notice that one always has the inequality $Q(u_1, u_2) \geq \frac{1}{4} |\nabla u_1|^2$, and hence 
$\dis{\eqref{eq:jw}}$ can be thought of as an extension of \eqref{eq:mt}. Our main result is the 
following one, which extends both \eqref{eq:mtsing} and \eqref{eq:jw}. 

\begin{thm}\label{t:main} 
Suppose $p_1, \dots, p_m \in \Sigma$ and  $\alpha_{i,j}$, $i = 1, 2$, $j = 1, \dots, m$, satisfy $\alpha_{i,j} 
> - 1$ for all $i, j$. Then, if $\wtilde{h}_i$ satisfy \eqref{eq:hi},  the following inequality holds 
\be \label{eq:mainineq} 4 \pi\sum_{i=1}^2 \min \left\{ 1, 1 + \min_j \alpha_{i,j} \right\} \log\int_{\Sigma} \wtilde{h}_i e^{u_i-\overline{u}_i}\,dV_{g} 
 \leq \int_{\Sigma} Q(u_1, u_2) \, dV_{g}+C \qquad \quad u_1, u_2 \in H^1(\Sigma).  
\ee
The constants in the above inequality are sharp. 
\end{thm}

\

\noindent  We expect the above result to be a main step for a possible variational approach for the study of \eqref{eq:e-1}. 
In the recent paper \cite{bjmr} the case of non-negative coefficients and positive genus has been treated using simply inequality \eqref{eq:jw}, 
as the  corresponding functions $\wtilde{h}_i$ are uniformly bounded (see also \cite{cheikh} and \cite{malruicpam} for the 
regular case). In more general cases, the full strength of \eqref{eq:mainineq} 
would be needed. 

\

\noindent  Some steps in the proof of the above theorem follow closely the arguments in \cite{jw}: through blow-up analysis one can show with few difficulties that  inequality $\dis{\eqref{eq:mainineq}}$ holds for any smaller couple of parameters, and moreover that there exist extremal functions for the corresponding 
Euler functionals \eqref{funzionale}. We pass then to the limit  for these extremals when the parameters approach the critical ones. 

However the presence of singularities might cause in principle a variety of blow-up behaviours (different blow-up rates for the two components, and 
blow-up at regular or singular points): using a Pohozaev identity from the recent paper \cite{linweizhang} we reduce ourselves 
to two cases only.  The former can be brought back to the scalar case, where one can use $\dis{\eqref{eq:mtsing}}$ to get a conclusion; the latter can be 
solved by using a \emph{local} version of the singular Moser-Trudinger inequality from Adimurthi and Sandeep \cite{as}. The latter argument 
in particular differs substantially from that in \cite{jw}, and it also provides a simpler argument for the regular case.

\

\section{Notation and preliminaries}

\noindent In this section we provide some useful notation and some known preliminary results which will be used in the proof of the main theorem.

First of all, given two points $\dis{x,y\in\Sg}$, we will indicate as $\dis{d(x,y)}$ the metric distance between $\dis{x}$ and $\dis{y}$ on $\dis{\Sg}$; we will denote as $\dis{B_r(p)}$ the open metric ball of radius $\dis{r}$ centered at $\dis{p}$. 

Given a function $\dis{u\in L^1(\Sg)}$, $\dis{\ov u}$ will stand for the average of $\dis{u}$ on $\dis{\Sg}$; since we will suppose, from now on, $\dis{|\Sg|=1}$, we can write
$$\ov u=\int_\Sg ud \, V_g.$$
We denote as $\dis{x^-}$ the negative part of a real number $\dis{x}$, that is $\dis{x^-:=\left\{\begin{array}{ll}0&\text{if }x\ge0\\-x&\text{if }x\le0\end{array}\right.}$, and we set, for $\dis{i\in\{1,2\}}$,
\be\label{eq:alfa}
\wtilde\a_i=-\max_{j\in\{1,\dots,m\}}{\a_{i,j}}^-.
\ee

Notice that, in these terms, the inequality we wish to prove is
$$4\pi\sum_{i=1}^2(1+\wtilde\a_i)\log\int_\Sg\wtilde h_ie^{u_i-\ov{u_i}}dV_g\le\int_\Sg Q(u_1,u_2)dV_g+C\qquad\quad u_1,u_2\in H^1(\Sg),$$
whereas the singular Chen-Troyanov $\dis{\eqref{eq:mtsing}}$ inequality can be expressed as
$$16\pi(1+\wtilde\a_i)\log\int_\Sg\wtilde h_ie^{u_i-\ov u_i}dV_g\le\int_\Sg\left|\nabla u\right|^2dV_g+C;\qquad\quad u\in H^1(\Sg).$$

We then define the $\dis{i^{th}}$ weight of a point $\dis{p\in\Sg}$, for $\dis{i\in\{1,2\}}$ in the following way
\be\label{eq:alfai}
p=q_j\quad\Rightarrow\quad \a_i(p)=\a_{i,j}\quad\quad p\notin\{q_1,\dots,q_m\}\quad\Rightarrow\quad\a_i(p)=0.
\ee
The definition implies that $\dis{\wtilde h_i\simeq d(\cdot,p)^{2\a_i(p)}}$ near $\dis{p}$; precisely, it is the only real number such that $\dis{\log\wtilde h_i-2\wtilde\a_i \log d(\cdot,p)}$ is bounded in a sufficiently small neighborhood of $\dis{p}$.\\
As anticipated in the introduction, we will prove inequality $\dis{\eqref{eq:mainineq}}$ via blow-up analysis. 
We define, for a sequence $\dis{u_k=(u_{1,k},u_{2,k})}$ of solutions of $\dis{\eqref{eq:todaregul}}$, the concentration value of the $\dis{i^{\text{th}}}$ component around a point $\dis{p\in\Sg}$ as
\be\label{eq:sigma}
\s_i(p):=\lim_{r\to0}\lim_{k\to+\infty}\int_{B_r(p)}\wtilde h_ie^{u_{i,k}}dV_g.
\ee
Lin, Wei and Zhang in \cite{linweizhang} found out, through a Poho\v zaev identity, that the concentration values satisfy the following condition, which was already pointed out for the regular case in \cite{jlw}.

\begin{thm}(\cite{linweizhang}, Proposition $\dis{3.1}$)
Let $\dis{u_k=(u_{1,k},u_{2,k})\in H^1(\Sg)^2}$ be solutions of $\dis{\eqref{eq:todaregul}}$, $\dis{\wtilde\a_i}$ be as in $\dis{\eqref{eq:alfai}}$ and $\dis{\s_i}$ be as in $\dis{\eqref{eq:sigma}}$. Then, it holds
\be\label{eq:poh}
\s_1(p)^2-\s_1(p)\s_2(p)+\s_2(p)^2=4\pi(1+\wtilde \a_1(p))\s_1(p)+4\pi(1+\wtilde \a_2(p))\s_2(p).
\ee
\end{thm}

\

\noindent In the setting we are considering, a dichotomy between concentration and compactness occurs, similar to the ones in the regular case from Jost-Wang \cite{jw}, Theorem $\dis{3.1}$. Since the proof of the theorem we are giving is very close to \cite{jw}, we will only sketch it; we refer to these papers for the details in the regular case.

\begin{thm}\label{t:comp} 	
Let $\dis{\wtilde h_i}$ as in $\dis{\eqref{eq:hi}}$, let $\dis{u_k=(u_{1,k},u_{2,k})\in H^1(\Sg)^2}$ be solutions of
$$\left\{\begin{array}{l}-\D u_{i,k}=2V_{i,k}\wtilde h_ie^{u_{i,k}}-V_{3-i,k}\wtilde h_{3-i}e^{u_{3-i,k}}+\psi_{i,k}\\\int_\Sg\wtilde h_ie^{u_{i,k}}dV_g\le C\\\|\psi_{i,k}\|_{L^p(\Sg)}\le C\\V_{i,k}\underset{k\to+\infty}\to1\text{ in }L^\infty(\Sg)\end{array}\right.\quad i\in\{1,2\},$$
for some $\dis{p>1}$, $\dis{C>0}$ and define the sets $\dis{S_i}$ as
$$S_i:=\left\{p\in\Sg:\;\exists\; x_k\underset{k\to+\infty}\to p\text{ such that }u_{i,k}(x_k)\underset{k\to+\infty}\to+\infty\right\}.$$
Then, after taking subsequences, one of the following alternatives happens.
\begin{enumerate}
\item For each $\dis{i\in\{1,2\}}$, either $\dis{u_{i,k}}$ is bounded in $\dis{L^\infty(\Sg)}$ or it tends uniformly to $\dis{-\infty}$.
\item $\dis{S_i\ne\emptyset}$ for some $\dis{i\in\{1,2\}}$; in this case, $\dis{S_i}$ is finite and either $\dis{u_{j,k}}$ is bounded in $\dis{L^\infty_{loc}(\Sg\backslash(S_1\cup S_2))}$ or it converges to $\dis{-\infty}$ in $\dis{L^\infty_{loc}(\Sg\backslash(S_1\cup S_2))}$ for each $\dis{j\in\{1,2\}}$; moreover, if $\dis{S_i\backslash S_{3-i}\ne\emptyset}$, then the latter alternative occurs for $\dis{u_{i,k}}$.
\end{enumerate}
\end{thm}

\begin{pfn} {\sc (Sketch)}
Reasoning as in \cite{bm} we find that, given $\dis{p\in\Sg}$, if for some $\dis{i\in\{1,2\}}$ one has 
$$
\dis{\limsup_{k\to+\infty}\int_{B_r(p)}V_{i,k}\wtilde h_ie^{u_{i,k}}dV_g<2\pi(1+\a_i(p)^-)}
$$ 
for sufficiently small $\dis{r}$, then $\dis{u_{i,k}}$ is uniformly bounded from above, and this fact implies the finiteness of the sets $\dis{S_i}$. The alternative between being bounded in $\dis{L^\infty}$ and converging uniformly to $\dis{-\infty}$ follows by applying a Harnack inequality and the last part of $\dis{(2)}$ follows by arguing as in \cite{breme}, Theorem $\dis{3}$.
\end{pfn}

\

\noindent Finally, as anticipated, we will need a singular Moser-Trudinger  inequality on bounded Euclidean domains, from \cite{as}:

\begin{thm}(\cite{as}, Theorem $\dis{2.1}$)
Let $\dis{\Omega\subset\R^2}$ a bounded domain containing the origin.
Then, for any $\dis{\a\in(-1,0]}$, it holds
$$\sup_{u\in H^1_0(\O),\int_\O|\n u(x)|^2dx\le1}\int_\O|x|^{2\a}e^{4\pi(1+\a)u(x)^2}dx\le C,$$
where $\dis{C}$ is a constant depending on $\dis{\a}$ and $\dis{\O}$ only.
\end{thm}

\

\noindent  From elementary inequalities we then obtain the following result. 

\begin{cor}
Let $\dis{\Omega\subset\R^2}$ a bounded domain containing the origin.
Then, for any $\dis{\a\in(-1,0]}$ and $\dis{u\in H^1_0(\O)}$, it holds
\be\label{eq:mte}
16\pi(1+\a)\log\int_\O|x|^{2\a}e^{u(x)}dx\le\int_\O|\n u(x)|^2dx+C.
\ee
\end{cor}

\

\section{A Moser-Trudinger inequality}

\noindent In this section, we are going to prove the following Moser-Trudinger type inequality.

\begin{thm}\label{t:mt<}
Let $\dis{\Sg}$ be a closed surface with area $\dis{|\Sg|=1}$, $\dis{\wtilde h_i}$ be as in $\dis{\eqref{eq:hi}}$, and $\dis{\wtilde\a_i}$ be as in $\dis{\eqref{eq:alfa}}$. 
Then, for any $\dis{\rh=(\rh_1,\rh_2)\in\R_+^2}$ satisfying $\dis{\rh_i<4\pi(1+\wtilde\a_i)}$ for both $\dis{i\in\{1,2\}}$ there exists $\dis{C(\rh)>0}$ such that the Euler-Lagrange functional $\dis{\eqref{funzionale}}$ verifies
$$J_\rh(u)>-C(\rh)\quad\quad\forall\;u\in H^1(\Sg)^2$$
\end{thm}

\begin{df}
As in \cite{jw}, we define the set of admissible parameters $\dis{\L}$ as
$$\L:=\left\{\rh\in\R_+^2:J_\rh\text{ is bounded from below}\right\}.$$
Clearly, $\dis{\L}$ preserves the partial order of $\dis{\R_+^2}$, that is if $\dis{\rh\in\L}$ then $\dis{\wtilde\rh\in\L}$ until  $\dis{\wtilde\rh_i\le\rh_i}$ for both $\dis{i\in\{1,2\}}$; in these terms, Theorem $\dis{\ref{t:mt<}}$ is equivalent to saying
$$(0,4\pi(1+\wtilde\a_1))\times(0,4\pi(1+\wtilde\a_2))\subset\L.$$
\end{df}

\begin{rem}
One can easily see that $\dis{\L}$ is not empty: since it holds
$$\frac{|\n u_1|^2+|\n u_2|^2}6\le Q(u_1,u_2)$$
one can apply the scalar Moser-Trudinger inequality $\dis{\eqref{eq:mtsing}}$ to both components to get
$$\left(0,\frac{8}3\pi(1+\wtilde\a_1)\right)\times\left(0,\frac{8}3\pi(1+\wtilde\a_2)\right)\subset\L.$$
\end{rem}

\noindent To prove Theorem $\dis{\ref{t:mt<}}$, some lemmas will be needed. First of all, we notice that when the parameter $\dis{\rh}$ is in the interior of the set $\dis{\L}$, then the energy functional is not only bounded from below, but even coercive and it has a minimizer; on the other hand, if $\dis{\rh}$ is on the boundary of $\dis{\L}$, then $\dis{J_\rh}$ cannot be coercive.

\begin{lem}\label{l:coe1}
For any $\dis{\rh\in\stackrel{\circ}\L}$ there exists a constant $\dis{C}$ such that
$$J_\rh(u)\ge\frac{\int_\Sg\left(|\n u_1|^2+|\n u_2|^2\right)dV_g}C-C$$
Moreover, $\dis{J_\rh}$ admits a minimizer $\dis{u=(u_1,u_2)}$ that solves $\dis{\eqref{eq:todaregul}}$.
\end{lem}

\begin{pf}
Taking $\dis{\d\in\left(0,\frac{d(\rh,\pa\L)}{\sqrt2}\right)}$, we have $\dis{(1+\d)\rh\in\L}$ so $\dis{J_{(1+\d)\rh}(u)\ge-C}$; therefore, we can write
\begin{eqnarray*}
J_\rh(u)&=&\frac{\d}{1+\d}\int_\Sg Q(u_1,u_2)dV_g+\frac{J_{(1+\d)\rh}(u)}{1+\d}\ge\\
&\ge&\frac{\d}{6(1+\d)}\int_\Sg\left(|\n u_1|^2+|\n u_2|^2\right)dV_g-C
\end{eqnarray*}
and the first claim follows.\\
To prove the rest we notice that, if we restrict ourselves to the subset of $\dis{H^1(\Sg)^2}$ consisting of all functions satisfying $\dis{\int_\Sg\wtilde h_ie^{u_i}dV_g=1}$, the energy is coercive because, from Poincar\'e's inequality and $\dis{\eqref{eq:mtsing}}$
\begin{eqnarray*}
\int_\Sg u_i^2dV_g&=&\int_\Sg\left(u_i-\ov{u_i}\right)^2dV_g+\left(\ov{u_i}\right)^2\le\\
&\le&C\int_\Sg|\n u_i|^2dV_g+\left(C+\frac{1}{16\pi(1+\wtilde\a_i)}\int_\Sg|\n u_i|^2dV_g\right)^2\le\\
&\le&C\left(1+\int_\Sg|\n u_i|^2dV_g\right)^2.
\end{eqnarray*}
Being $\dis{J_\rh}$ weakly lower-semicontinuous as well, the existence of minimizers follows from the direct methods of calculus of variations.
\end{pf}

\begin{lem}\label{l:coe2}
For any $\dis{\rh\in\pa\L}$ there exists a sequence $\dis{\{\wtilde u_k\}_{k\in\N}\subset H^1(\Sg)^2}$ verifying
$$\int_\Sg\left(|\n\wtilde u_{1,k}|^2+|\n\wtilde u_{2,k}|^2\right)dV_g\underset{k\to+\infty}\to+\infty\quad\quad\lim_{k\to+\infty}\frac{J_\rh(\wtilde u_k)}{\int_\Sg\left(|\n\wtilde u_{1,k}|^2+|\n\wtilde u_{2,k}|^2\right)dV_g}\le0.$$
\end{lem}

\begin{pf}
Suppose by contradiction that
$$\int_\Sg\left(|\n u_{1,k}|^2+|\n u_{2,k}|^2\right)dV_g\underset{k\to+\infty}\to+\infty\quad\quad\Rightarrow\quad\quad\frac{J_\rh(u_k)}{\int_\Sg\left(|\n u_{1,k}|^2+|\n u_{2,k}|^2\right)dV_g}\ge\th>0$$
for any choice of $\dis{\{u_k\}}$. This would mean that $\dis{J_\rh(u)\ge\frac{\th}2\int_\Sg\left(|\n u_1|^2+|\n u_2|^2\right)dV_g-C}$, hence for any small $\dis{\d}$ we would get
\begin{eqnarray*}
J_{(1+\d)\rh}(u)&=&(1+\d)J_\rh(u)-\d\int_\Sg Q(u_1,u_2)dV_g\ge\\
&\ge&\left((1+\d)\frac{\th}2-\frac{\d}2\right)\int_\Sg\left(|\n u_1|^2+|\n u_2|^2\right)dV_g-C\\
&\ge&-C
\end{eqnarray*}
hence $\dis{(1+\d)\rh\in\L}$, whereas one clearly has $\dis{(1-\d)\rh\in\L}$; this is in contradiction to $\dis{\rh\in\pa\L}$.
\end{pf}

\

\noindent We then need a basic calculus lemma. Its proof will be omitted, as it can be found in \cite{jw} (following an idea of W. Ding).

\begin{lem}[\cite{jw}, Lemma $\dis{4.4}$]\label{l:fakbk}
Let $\dis{\{a_k\}_{k\in\N}}$ and $\dis{\{b_k\}_{k\in\N}}$ be two sequences of real numbers satisfying
$$a_k\underset{k\to+\infty}\to+\infty\quad\quad\text{and}\quad\quad\lim_{k\to+\infty}\frac{b_k}{a_k}\le0.$$
Then there exists a smooth function $\dis{F:[0,+\infty)\to\R}$ satisfying, up to subsequences,
$$0<F'(t)<1\quad\text{for any }t\ge0\quad\quad F'(t)\underset{t \to+\infty}\to0\quad\quad F(a_k)-b_k\underset{k\to+\infty}\to+\infty.$$
\end{lem}

\noindent The latter lemma will be applied to the sequences
$$a_k=\int_\Sg Q(\wtilde u_{1,k},\wtilde u_{2,k})dV_g\quad\quad\quad\quad b_k=J_\rh(\wtilde u_k)$$
where $\dis{\wtilde u_k}$ is as in Lemma $\dis{\ref{l:coe2}}$, and we will consider the auxiliary functional
$$\wtilde J_\rh(u):=J_\rh(u)-F\left(\int_\Sg Q(u_1,u_2)dV_g\right),$$
whose behavior is described by the following lemma.

\begin{lem}\label{l:irho}
For any $\dis{\rh\in\stackrel{\circ}\L}$ the functional $\dis{\wtilde J_\rh}$ is bounded from below on $\dis{H^1(\Sg)^2}$ and its infimum is achieved by a function satisfying
$$\left\{\begin{array}{l}-\left(1-\frac{2}3g(u)\right)\D u_i+\frac{g(u)}3\D u_{3-i}=2\rh_i\left(\wtilde h_ie^{u_i}-1\right)-\rh_{3-i}\left(\wtilde h_{3-i}e^{u_{3-i}}-1\right); \\\int_\Sg\wtilde h_ie^{u_i}dV_g=1,\end{array}\right.$$
where $\dis{g(u)=F'\left(\int_\Sg Q(u_1,u_2)dV_g\right)}$.
On the other hand, if $\dis{\rh\in\pa\L}$ then $\dis{\inf_{H^1(\Sg)^2}\wtilde J_\rh=-\infty}$
\end{lem}

\begin{pf}
For $\dis{\rh\in\stackrel{\circ}\L}$ one can argue as in Lemma $\dis{\ref{l:coe1}}$, yielding lower semi-continuity from the regularity of $\dis{F}$ and coerciveness from the behavior of $\dis{F'}$ at infinity.\\
For $\dis{\rh\in\pa\L}$, taking $\dis{\wtilde u_k}$ as in Lemma $\dis{\ref{l:coe2}}$ and applying Lemma $\dis{\ref{l:fakbk}}$ one gets
$$\wtilde J_\rh(\wtilde u_k)=b_k-F(a_k)\underset{k\to+\infty}\to-\infty.$$
This concludes the proof.
\end{pf}

\

\noindent We are now in position to prove the main theorem of this section.

\

\begin{pfn}\begin{sc} of Theorem $\dis{\ref{t:mt<}}$\end{sc}
Suppose by contradiction that
$$(0,4\pi(1+\wtilde\a_1))\times(0,4\pi(1+\wtilde\a_2))\not\subset\L;$$
then there is some $\dis{\ov\rh\in\pa\L}$ with $\dis{\ov\rh_i<4\pi(1+\wtilde\a_i)}$ for both $\dis{i\in\{1,2\}}$.\\
Consider a sequence $\dis{\{\rh_k\}_{k\in\N}\in\stackrel{\circ}\L}$ with $\dis{\rh_k\underset{k\to+\infty}\to\ov\rh}$ and a minimizer $\dis{u_k}$ for $\dis{I_{\rh_k}}$, as in Lemma $\dis{\ref{l:irho}}$; then, $\dis{v_k:=u_k+\log\rh_k}$ solves
$$\left\{\begin{array}{l}-\D v_{i,k}=2\frac{6-5g(v_k)}{6-8g(v_k)+2g(v_k)^2}\left(\wtilde h_ie^{v_{i,k}}-\rh_{i,k}\right)-\frac{3-4g(v_k)}{3-4g(v_k)+g(v_k)^2}\left(\wtilde h_{3-i}e^{v_{3-i,k}}-\rh_{3-i,k}\right); \\\int_\Sg\wtilde h_ie^{v_{i,k}}dV_g=\rh_{i,k},\end{array}\right.$$
with $\dis{\frac{6-5g(v_k)}{6-8g(v_k)+2g(v_k)^2}}$ and $\dis{\frac{3-4g(v_k)}{3-4g(v_k)+g(v_k)^2}}$ both uniformly converging to $\dis{1}$, so Theorem $\dis{\ref{t:comp}}$ can be applied to this sequence. The normalization on the integral implies that $\dis{u_{i,k}}$ cannot tend to $\dis{-\infty}$ for any $\dis{i\in\{1,2\}}$; moreover, we can also exclude boundedness in $\dis{L^\infty(\Sg)}$ because this would imply convergence to a minimizer $\dis{\ov u}$ of $\dis{I_{\ov\rh}}$, contradicting Lemma $\dis{\ref{l:irho}}$.\\

The only case left is the blow-up around at least one point $\dis{p}$: Poho\v zaev's identity $\dis{\eqref{eq:poh}}$ implies that if there is a singularity of mass $\dis{\a_{i,j}}$ on $\dis{p}$ then $\dis{\sigma_i\ge4\pi(1+\alpha_{i,j})}$ for some $\dis{i\in\{1,2\}}$, whereas if $\dis{p}$ is a regular point then there is a component with a mass of at least $\dis{4\pi}$ around it; in both cases, for such an $\dis{i}$ we obtain:
$$4\pi(1+\wtilde\a_i)\le\lim_{r\to0}\lim_{k\to+\infty}\int_{B_r(p)}\wtilde h_ie^{v_{i,k}}dV_g\le\lim_{k\to+\infty}\int_\Sg\wtilde h_ie^{v_{i,k}}dV_g=\ov\rh_i<4\pi(1+\wtilde\a_i),$$
that is a contradiction.
\end{pfn}

We conclude the section by showing a partial converse of Theorem $\dis{\ref{t:mt<}}$, namely that for higher values of the parameter $\dis{\rh}$ the functional $\dis{J_\rh}$ is unbounded from below.

\begin{pro}\label{p:mt>}
If $\dis{\rh_i>4\pi(1+\wtilde\a_i)}$ for some $\dis{i\in\{1,2\}}$, then $\dis{\inf_{H^1(\Sg)^2}J_\rh=-\infty}$ that is
$$\L\subset(0,4\pi(1+\wtilde\a_1)]\times\left(0,4\pi\left(1+\wtilde\a_2\right)\right].$$
\end{pro}

\begin{pf}
We will show the proof only for $\dis{i=1}$, since the same argument works for $\dis{i=2}$ as well.\\
Choosing a point $\dis{p_1}$ such that $\dis{\wtilde h_1\simeq d(\cdot,p_i)^{2\wtilde\a_1}}$ in its neighborhood, we define for large $\dis{\l}$
$$\var_{1,\l}(x)=\log\left(\frac{\l^{1+\wtilde\a_1}}{1+(\l d(x,p_1))^{2(1+\wtilde\a_1)}}\right)^2; \quad\quad \var_{2,\l}(x)=-\frac{1}2\log\left(\frac{\l^{1+\wtilde\a_1}}{1+(\l d(x,p_1))^{2(1+\wtilde\a_1)}}\right)^2. $$
Using the fact that $\dis{\left|\n\left(d(x,p_1)^{2(1+\wtilde\a_1)}\right)\right|\le2(1+\wtilde\a_1)d(x,p_1)^{1+2\wtilde\a_1}}$, we obtain
\begin{eqnarray*}
|\n\var_{1,\l}(x)|&=&\left|\frac{-2\l^{2(1+\wtilde\a_1)}\left|\n\left(d(x,p_1)^{2(1+\wtilde\a_1)}\right)\right|}{1+(\l d(x,p_1))^{2(1+\wtilde\a_1)}}\right|\le\\
&\le&\frac{4(1+\wtilde\a_1)\l^{2(1+\wtilde\a_1)}d(x,p_1)^{1+2\wtilde\a_1}}{1+(\l d(x,p_1))^{2(1+\wtilde\a_1)}}\le\\
&\le&\min\left\{C\l,\frac{4(1+\wtilde\a_1)}{d(x,p_1)}\right\},
\end{eqnarray*}
and therefore
\begin{eqnarray}
\nonumber\int_\Sg Q(\var_{1,\l},\var_{2,\l})dV_g&=&\frac{1}4\int_\Sg|\n\var_{1,\l}|^2dV_g\le\\
&\le&C\l^2\int_{B_{\frac{1}\l}(p_1)}dV_g+4(1+\wtilde\a_1)^2\int_{\Sg\backslash B_{\frac{1}\l}(p_1)}\frac{dV_g}{d(\cdot,p_1)^2}\le\\
\label{eq:grad}&\le&C+8\pi(1+\wtilde\a_1)^2\log\l. \nonumber 
\end{eqnarray}
Moreover, being
\be\label{eq:stima}
\max\{1,(\l d(x,p_1))^{2(1+\wtilde\a_1)}\}\le 1+(\l d(x,p_1))^{2(1+\wtilde\a_1)}\le C\max\{1,(\l d(x,p_1))^{2(1+\wtilde\a_1)}\},
\ee
one gets
$$\nonumber\ov{\var_{1,\l}}=\int_\Sg(\max\{2(1+\wtilde\a_1)\log\l,-2(1+\wtilde\a_1)(\log\l+2\log d(\cdot,p_1))\}+O(1))dV_g.$$
Dividing $\dis{\Sg}$ into the two regions where the above maximum is attained and using the integrability of $\dis{\log d(\cdot,p_1)}$ in two dimensions one gets
\begin{eqnarray}
\nonumber\ov{\var_{1,\l}}&=&2(1+\wtilde\a_1)\log\l\int_{B_\frac{1}\l(p_1)}dV_g-2(1+\wtilde\a_1)\log\l\int_{\Sg\backslash B_\frac{1}\l(p_1)}dV_g-\\
\label{eq:media}&-&4(1+\wtilde\a_1)\int_{\Sg\backslash B_\frac{1}\l(p_1)}\log d(\cdot,p_1)dV_g+O(1)\\
\nonumber&=&-2(1+\wtilde\a_1)\log\l+O(1),
\end{eqnarray}
and clearly $\dis{\ov{\var_{2,\l}}=(1+\wtilde\a_1)\log\l+O(1)}$.

For a small but fixed $\d > 0$ we have, again by $\dis{\eqref{eq:stima}}$,
\begin{eqnarray}
\nonumber\int_\Sg\wtilde h_1e^{\var_{1,\l}}dV_g&\ge&C\int_{B_\d(p_1) \backslash B_{\frac{1}\l}(p_1)}d(\cdot,p_1)^{2\wtilde\a_1}e^{\var_{1,\l}}dV_g\ge\\
&\ge&\frac{C}{\l^{2(1+\wtilde\a_1)}}\int_{B_\d(p_1) \backslash B_\frac{1}\l(p_1)}\frac{dV_g}{d(\cdot,p_1)^{4+2\wtilde\a_1}}\ge\\
\label{eq:exp1}&\ge&C; \nonumber 
\end{eqnarray}
on the other hand, we can write that
\begin{eqnarray}
\nonumber\int_\Sg\wtilde h_2e^{\var_{2,\l}}dV_g&\ge&C\l^{1+\wtilde\a_1}\int_{\Sg\backslash B_\frac{1}\l(p_1)}\wtilde h_2d(\cdot,p_1)^{2(1+\wtilde\a_1)}dV_g\ge\\
\label{eq:exp2}&\ge&C\l^{1+\wtilde\a_1}.
\end{eqnarray}
Therefore, from $\dis{\eqref{eq:grad},\eqref{eq:media},\eqref{eq:exp1},\eqref{eq:exp2}}$ we conclude that
$$J_\rh(u)\le2(1+\wtilde\a_1)(4\pi(1+\wtilde\a_1)-\rh_1)\log\l+O(1)\underset{\l\to\infty}\to-\infty,$$
as desired.
\end{pf}

\

\section{The optimal inequality}

\noindent In the last section we are going to discuss the boundedness from below of $\dis{J_\rho}$ in the only case left, that is when $\dis{\rh_i=4\pi(1+\wtilde\a_i)}$ for some $\dis{i\in\{1,2\}}$; we will show that $\dis{\inf_{H_1(\Sg)^2}J_\rh>-\infty}$ in this case as well.

\begin{thm}\label{t:mt=}
Let $\dis{\Sg}$ be a closed surface with area $\dis{|\Sg|=1}$, $\dis{\wtilde h_i}$ be as in $\dis{\eqref{eq:hi}}$, $\dis{\wtilde\a_i}$ be as in $\dis{\eqref{eq:alfa}}$ and $\dis{J_\rho}$ be as in $\dis{\eqref{funzionale}}$.\\
Then, there exists a constant $\dis{C>0}$ such that for any $\dis{u\in H^1(\Sg)^2}$
$$J_{4\pi(1+\wtilde\a_1),4\pi(1+\wtilde\a_2)}(u)>-C$$
namely
$$\L=(0,4\pi(1+\wtilde\a_1)]\times(0,4\pi(1+\wtilde\a_2)].$$
\end{thm}

\

\noindent Theorem $\dis{\ref{t:mt=}}$ is equivalent to saying that, given a sequence $\dis{\rh_k\underset{k\to+\infty}\nearrow(4\pi(1+\wtilde\a_1),4\pi(1+\wtilde\a_2))}$, there exists $\dis{C>0}$ such that $\dis{\inf_{H^1(\Sg)^2}J_{\rh_k}\ge-C}$.\\
Moreover, in view of Lemma $\dis{\ref{l:coe1}}$, it suffices to show that the minimizers $\dis{u_k}$ of $\dis{J_{\rh_k}}$ verify $\dis{J_{\rh_k}(u_k)>-C}$; these functions solve
$$\left\{\begin{array}{l}-\D u_{i,k}=2\rh_{i,k}\left(\wtilde h_ie^{u_{i,k}}-1\right)-\rh_{3-{i,k}}\left(\wtilde h_{3-i}e^{u_{3-i,k}}-1\right)\\\int_\Sg\wtilde h_ie^{u_{i,k}}dV_g=1\end{array}\right.\quad i\in\{1,2\},$$
therefore, as in the proof of Theorem $\dis{\ref{t:mt<}}$, we can apply Theorem $\dis{\ref{t:comp}}$ to $\dis{v_k:=u_k+\log\rh_k}$.

As in the proof of Theorem \ref{t:mt<}, the condition on the integral excludes convergence to $\dis{-\infty}$, whereas if $\dis{u_k}$ is bounded in $\dis{\|\cdot\|_{L^\infty(\Sg)}}$ it converges to a minimizer of $\dis{J_{4\pi(1+\wtilde\a_1),4\pi(1+\wtilde\a_2)}}$ hence the conclusion is trivial, so we may suppose that at least one component blows up.\\
The following lemma describes the two possible blow-up scenarios.

\begin{lem}\label{l:blow-up}
Let $\dis{\{u_k\}_{k\in\N}}$ be a blowing up sequence of minimizers of $\dis{J_{\rh_k}}$ for some sequence $\rho_k$ such that  $\dis{\rh_k\underset{k\to+\infty}\to(4\pi(1+\wtilde\a_1),4\pi(1+\wtilde\a_2))}$ and let $\dis{\a_i(p)}$ be as in $\dis{\eqref{eq:alfai}}$. Then, one of the following happens:
\begin{enumerate}
\item Only the $\dis{i^{\text{th}}}$ component of $\dis{u_k}$ blows up, for some $\dis{i\in\{1,2\}}$ and it does at a single point $\dis{p_i}$ with $\dis{\a_i(p_i)=\wtilde\a_i}$ around it.
\item Each component of $\dis{u_k}$ blows up at a single point $\dis{p_i}$ satisfying $\dis{\a_i(p_i)=\wtilde\a_i}$ around it, and $\dis{p_1\ne p_2}$.
\end{enumerate}
\end{lem}
\begin{pf}
Suppose that only one component blows up, say $\dis{u_{1,k}}$, and suppose it blows up around a point $\dis{p_1}$ satisfying $\dis{\a_1(p_1)>\wtilde\a_1}$. Then, by $\dis{\eqref{eq:poh}}$ we obtain
$$
4\pi(1+\a_1(p_1))=\lim_{r\to0}\lim_{k\to+\infty}\int_{B_r(p_1)}\wtilde h_1 e^{v_{1,k}}dV_g\le\lim_{k\to+\infty}\int_\Sg\wtilde h_1 e^{v_{1,k}}dV_g=4\pi(1+\wtilde\a_1),
$$
that is a contradiction; moreover, if the blow-up occurs at two points $\dis{p_1}$, $\dis{\ov p_1}$, then one similarly gets another contradiction:
$$
8\pi(1+\wtilde\a_1)=\lim_{r\to0}\lim_{k\to+\infty}\int_{B_r(p_1)\cup B_r(p_2)}\wtilde h_1 e^{v_{1,k}}dV_g\le\lim_{k\to+\infty}\int_\Sg\wtilde h_1 e^{v_{1,k}}dV_g=4\pi(1+\wtilde\a_1).
$$
Suppose now that both components blow up at the same point; then, again by $\dis{\eqref{eq:poh}}$, $\dis{v_{i,k}}$ must have a local mass strictly greater than $\dis{4\pi(1+\wtilde\a_i)}$ around that point, for some $\dis{i\in\{1,2\}}$, but this is impossible since the total mass of $\dis{v_{i,k}}$ is converging to $\dis{4\pi(1+\wtilde\a_i)}$; therefore, at any given point only one component may blow up, hence we can argue as in the previous case to get the conclusion. 
\end{pf}

We will consider first the single-component blow-up in alternative $\dis{(1)}$.
\begin{lem}\label{l:g1}
Suppose $\dis{u_{1,k}}$ blows up at $\dis{p_1}$ and $\dis{u_{2,k}}$ does not blow up. Then,
\begin{enumerate}
\item $\dis{u_{1,k}-\overline{u_{1,k}}\underset{k\to+\infty}\to G_1}$ in $\dis{W^{2,p}_{loc}(\Sg\backslash\{p_1\})}$ for any $\dis{p\in\left[1,\frac{1}{-\wtilde\a_1}\right)}$ and weakly$\dis{^*}$ in $\dis{W^{1,q}(\Sg)}$ for any $\dis{q\in[1,2)}$, and $\dis{G_1}$ satisfies
\be\label{eq:g1}
\left\{\begin{array}{l}-\D G_1=8\pi(1+\wtilde\a_1)\left(\d_{p_1}-1\right)-4\pi(1+\wtilde\a_2)(f-1)\\\int_\Sg G_1dV_g=0\end{array}\right..
\ee
\item $\dis{u_{2,k}-\ov{u_{2,k}}\underset{k\to+\infty}\to G_2}$ in $\dis{W^{2,p}_{loc}(\Sg\backslash\{p_1\})}$ for any $\dis{p\in\left[1,\frac{1}{-\wtilde\a_2}\right)}$ and weakly$\dis{^*}$ in $\dis{W^{1,q}(\Sg)}$ for any $\dis{q\in[1,2)}$, and $\dis{G_2}$ satisfies
\be\label{eq:g2}
\left\{\begin{array}{l}-\D G_2=8\pi(1+\wtilde\a_2)(f-1)-4\pi(1+\wtilde\a_1)\left(\d_{p_1}-1\right)\\\int_\Sg G_2dV_g=0\end{array}\right..
\ee
for some non-negative $\dis{f\in L^1(\Sg)}$ satisfying $\dis{\int_\Sg fdV_g=1}$.
\end{enumerate}
\end{lem}

\begin{pf}
First of all, we prove that $\dis{u_{i,k}-\ov{u_{i,k}}}$ is bounded in $\dis{W^{1,q}(\Sg)}$ for $\dis{q\in[1,2)}$: taking $\dis{q'\in(2,+\infty]}$ such that $\dis{\frac{1}{q'}+\frac{1}q=1}$,
\begin{eqnarray*}
\|u_{i,k}-\ov{u_{i,k}}\|_{W^{1,q}(\Sg)}&\le&C\|\n u_{i,k}\|_{L^q(\Sg)}=\\
&=&C\sup_{\phi\in W^{1,q'}(\Sg),\;\|\n\phi\|_{L^{q'}}\le1}\left|\int_\Sg\n u_{i,k}\cdot\n\phi dV_g\right|\le\\
&\le&C\sup_{\phi\in W^{1,q'}(\Sg),\;\|\n\phi\|_{L^{q'}}\le1}\|\D u_{i,k}\|_{L^1(\Sg)}\|\phi\|_{L^\infty(\Sg)}\le\\
&\le&C\sup_{\phi\in W^{1,q'}(\Sg),\;\|\n\phi\|_{L^{q'}}\le1}\|\D u_{i,k}\|_{L^1(\Sg)}\|\n\phi\|_{L^{q'}(\Sg)}\le\\
&\le&C.
\end{eqnarray*}
Moreover, from Theorem $\dis{\ref{t:comp}}$ we know that, in the sense of measure,
$$\wtilde h_1e^{u_{1,k}}\underset{k\to+\infty}\wk\d_{p_1}\quad\quad\wtilde h_2e^{u_{2,k}}\underset{k\to+\infty}\wk f\in L^1(\Sg);$$
therefore, taking $\dis{G_i}$ satisfying respectively $\dis{\eqref{eq:g1}}$, ${\eqref{eq:g2}}$, for any fixed $\dis{\phi\in W^{1,q'}(\Sg)}$
\begin{eqnarray*}
\left|\int_\Sg\n\left(u_{1,k}-\ov{u_{1,k}}-G_1\right)\cdot\n\phi dV_g\right|&=&\int_\Sg\left(-\D u_{1,k}+\D G_1\right)\phi dV_g\le\\
&\le&C\left|\int_\Sg\left(2\rh_{1,k}\wtilde h_1e^{u_{1,k}}-8\pi(1+\wtilde\a_1)\d_{p_1}\right)\phi dV_g\right|+\\
&+&C\left|\int_\Sg\left(4\pi(1+\wtilde\a_2)f-\rh_{2,k}\wtilde h_2e^{u_{2,k}}\right)\phi dV_g\right|=\\
&=&o(1).
\end{eqnarray*}
in a similar way, we get $\dis{u_{2,k}-\ov{u_{2,k}}\underset{k\to+\infty}{^*\wk}G_2}$ in $\dis{W^{1,q}(\Sg)}$ and convergence in $\dis{W^{2,p}_{loc}(\Sg\backslash\{p_1\})}$ follows from standard elliptic estimates.
\end{pf}

\begin{rem}
From the previous lemma, we deduce that $\dis{|\ov{u_{2,k}}|\le C}$, since both $\dis{u_{2,k}}$ and $\dis{u_{2,k}-\ov{u_{2,k}}}$ are uniformly bounded in $\dis{L^\infty_{loc}(\Sg\backslash\{p_1\})}$; therefore, up to subsequences, the previous convergence result extends to $\dis{u_{2,k}}$.
\end{rem}

\

%The following lemma show that, in this case, a particular linear combination of $\dis{u_{1,k}}$ and $\dis{u_{2,k}}$ is negligible for $\dis{J_{\rh_k}}$.

%\begin{lem}\label{l:vk1}
%Suppose $\dis{u_{1,k}}$ blows up at $\dis{p_1}$ and $\dis{u_{2,k}}$ does not blow up. Then, the sequence
%$$v_{2,k}:=2(u_{2,k}-\ov{u_{2,k}})+u_{1,k}-\ov{u_{1,k}}$$
%is uniformly bounded in $\dis{W^{2,p}(\Sg)}$ for any $\dis{p\in\left[1,\frac{1}{-\wtilde\a_1}\right)}$.
%\end{lem}
%\begin{pf}
%$\dis{v_{2,k}}$ solves
%$$\left\{\begin{array}{l}-\D v_{2,k}=3\rh_{2,k}\left(\wtilde h_2e^{u_{2,k}}-1\right)\\\int_\Sg v_{2,k}=0\end{array}\right.$$
%therefore the thesis follows from the boundedness from above of $\dis{u_{2,k}}$ and standard elliptic estimates.
%\end{pf}

\noindent We will now consider the alternative $\dis{(2)}$ in Lemma $\dis{\eqref{l:blow-up}}$. 

When both components blow up, the last lemma have a counterpart; its proof follow closely the proof of Lemma $\dis{\ref{l:g1}}$, and therefore will be omitted.

\begin{lem}\label{l:g2}
Suppose each $\dis{u_{i,k}}$ blows up at $\dis{p_i}$. Then, for both $\dis{i\in\{1,2\}}$ we have that 
$\dis{u_{i,k}-\overline{u_{i,k}}\underset{k\to+\infty}\to G_i}$ in $\dis{W^{2,p}_{loc}(\Sg\backslash\{p_i\})}$ for any 
$\dis{p\in\left[1,\frac{1}{-\wtilde\a_i}\right)}$ ($p \in [1,\infty)$ if $\wtilde{\a}_i = 0$) and weakly$\dis{^*}$ in $\dis{W^{1,q}(\Sg)}$ for any $\dis{q\in[1,2)}$, and $\dis{G_i}$ satisfies
$$
\left\{\begin{array}{l}-\D G_i=8\pi(1+\wtilde\a_i)\left(\d_{p_i}-1\right)-4\pi(1+\wtilde\a_{3-i})(\d_{p_{3-i}}-1); \\\int_\Sg G_idV_g=0.\end{array}\right.
$$
\end{lem}

%\begin{lem}\label{l:vk2}
%Suppose each $\dis{u_{i,k}}$ blows up at $\dis{p_i}$. Then, for both $\dis{i\in\{1,2\}}$, the sequence
%$$v_{i,k}:=2(u_{i,k}-\ov{u_{i,k}})+u_{3-i,k}-\ov{u_{3-i,k}}$$
%is uniformly bounded in $\dis{W^{2,p}_{loc}(\Sg\backslash\{p_i\})}$ for any $\dis{p\in\left[1,\frac{1}{-\wtilde\a_i}\right)}$.
%\end{lem}

\

\noindent In the case of both components blowing up, a sort of \emph{localized} Moser-Trudinger inequality is required.

\begin{lem}\label{l:mtloc}
Suppose each $\dis{u_{i,k}}$ blows up at $\dis{p_i}$. Then, for any small $\dis{\d>0}$ there exists $\dis{C=C(\d)>0}$ such that for both $\dis{i\in\{1,2\}}$
$$\frac{1}4\int_{B_\d(p_i)}|\n u_{i,k}|^2dV_g+\rh_{i,k}\ov{u_{i,k}}\ge-C.$$
\end{lem}

\begin{pf}
We will take $\dis{\d}$ such that $\dis{B_\d(p_i)}$ does not contain any other singular point and we will suppose that $\dis{B_\d(p_i)}$ is a flat disk, see \cite{jw} (Remark $\dis{3.3}$).\\

This condition can be achieved through a conformal change of metric which results in a modified Liouville equation. The same estimates on minimizers hold true for the modified equation and one gets lower bounds on the functionals as before. 

Consider the solution $\dis{\wtilde w_{i,k}}$ of
$$\left\{\begin{array}{ll}-\D\wtilde w_{i,k}=0&\text{on }B_\d(p_i), \\\wtilde w_{i,k}-u_{i,k}+\ov{u_{i,k}}=0&\text{on }\pa B_\d(p_i); \end{array}\right.$$
standard elliptic estimates and Lemma $\dis{\ref{l:g2}}$ give
$$\left\|\wtilde w_{i,k}\right\|_{C^1(B_\d(p_i))}\le C\left\|\wtilde w_{i,k}\right\|_{L^\infty(B_\d(p_i))}\le C\|u_{i,k}-\overline{u_{i,k}}\|_{L^\infty(\pa B_\d(p_i))}\le C.$$
Moreover, we can apply the scalar Moser-Trudinger inequality $\dis{\eqref{eq:mte}}$ to $\dis{w_{i,k}:=u_{i,k}-\ov{u_{i,k}}-\wtilde w_{i,k}}$, which belongs to $\dis{H^1_0(B_\d(p_i))}$:
$$\int_{B_\d(p_i)}|\n w_{i,k}|^2dV_g-16\pi(1+\wtilde\a_i)\log\int_{B_\d(p_i)}d(\cdot,p_i)^{2\wtilde\a_i}e^{w_{i,k}}dV_g\ge-C.$$
The construction of $\dis{\wtilde w_{i,k}}$ gives
\begin{eqnarray*}
\int_{B_\d(p_i)}|\n w_{i,k}|^2dV_g-\int_{B_\d(p_i)}|\n u_{i,k}|^2dV_g&=&\int_{B_\d(p_i)}\left(2\n u_{i,k}\cdot\n\wtilde w_{i,k}+|\n\wtilde w_{i,k}|^2\right)dV_g\le\\
&\le&2|\n\wtilde w_{i,k}|_{L^\infty(B_\d(p_i))}\int_{B_\d(p_i)}|\n u_{i,k}|dV_g+\\
&+&\int_{B_\d(p_i)}|\n\wtilde w_{i,k}|^2dV_g\le\\
&\le&C;\\
\end{eqnarray*}
on the other hand, for large $\dis{k}$ we may suppose that $\dis{\int_{B_\d(p_i)}\wtilde h_ie^{u_{i,k}}dV_g\ge\frac{1}2}$, so
\begin{eqnarray*}
\int_{B_\d(p_i)}d(\cdot,p_i)^{2\wtilde\a_i}e^{w_{i,k}}dV_g&=&e^{-\overline{u_{i,k}}}\int_{B_\d(p_i)}d(\cdot,p_i)^{2\wtilde\a_i}e^{u_{i,k}-\wtilde w_{i,k}}dV_g\ge\\
&\ge&Ce^{-\overline{u_{i,k}}}\int_{B_\d(p_i)}\wtilde h_ie^{u_{i,k}-\wtilde w_{i,k}}dV_g\ge\\
&\ge&Ce^{-\overline{u_{i,k}}}\int_{B_\d(p_i)}\wtilde h_ie^{u_{i,k}}dV_g\ge\\
&\ge&\frac{C}2e^{-\overline{u_{i,k}}}.
\end{eqnarray*}
Therefore, we get 
\begin{eqnarray*}
\frac{1}4\int_{B_\d(p_i)}|\n u_{i,k}|^2dV_g+\rh_{i,k}\ov{u_{i,k}}&\ge&\frac{1}4\int_{B_\d(p_i)}|\n w_{i,k}|^2dV_g-\\
&-&\rh_{i,k} \log\int_{B_\d(p_i)}d(\cdot,p_i)^{2\wtilde\a_i}e^{w_{i,k}}dV_g-C\ge\\
&\ge&-C.ì
\end{eqnarray*}
which is the conclusion
\end{pf}

\

\noindent We have now all the necessary tools to conclude the proof of Theorem $\dis{\ref{t:mt=}}$.

\

\begin{pfn}\begin{sc} of Theorem $\dis{\ref{t:mt=}}$\end{sc}
Take a minimizing blowing up sequence $\dis{u_k}$ and suppose that the first alternative in Lemma $\dis{\ref{l:blow-up}}$ holds; it is not restrictive to suppose that $\dis{u_{1,k}}$ blows up.

From Lemma $\dis{\ref{l:g1}}$ and the following remark we know that $\dis{\ov{u_{2,k}}}$ is uniformly bounded; therefore, using the scalar Moser-Trudinger inequality $\dis{\eqref{eq:mtsing}}$ we obtain
\begin{eqnarray*}
J_{\rh_k}(u_k)&=&\int_\Sg Q(u_{1,k},u_{2,k})dV_g+\rh_{1,k}\ov{u_{1,k}}+\rh_{2,k}\ov{u_{2,k}}\ge\\
&\ge&\int_\Sg Q(u_{1,k},u_{2,k})dV_g+\rh_{1,k}\ov{u_{1,k}}-C\ge\\
&\ge&\frac{1}4\int_\Sg|\n u_{1,k}|^2dV_g+\rh_{1,k}\ov{u_{1,k}}-C\ge\\
&\ge&-C.
\end{eqnarray*}
that concludes the analysis of the first case.\\
Suppose now that both components blow up; then, we may conclude by applying Lemma $\dis{\ref{l:mtloc}}$:
\begin{eqnarray*}
J_{\rh_k}(u_k)&=&\int_\Sg Q(u_{1,k},u_{2,k})dV_g+\rh_{1,k}\ov{u_{1,k}}+\rh_{2,k}\ov{u_{2,k}}\ge\\
&\ge&\sum_{i=1}^2\left(\int_{B_\d(p_i)}Q(u_{1,k},u_{2,k})dV_g+\rh_{i,k}\ov{u_{i,k}}\right)\ge\\
&\ge&\sum_{i=1}^2\left(\frac{1}4\int_{B_\d(p_i)}|\n u_{i,k}|^2dV_g+\rh_{i,k}\ov{u_{i,k}}\right)-C\ge\\
&\ge&-C.
\end{eqnarray*}
This concludes the proof.
\end{pfn}

\bibliography{mttoda}

\begin{thebibliography}{10}

\bibitem{as}
Adimurthi and K.~Sandeep.
\newblock A singular {M}oser-{T}rudinger embedding and its applications.
\newblock {\em NoDEA Nonlinear Differential Equations Appl.}, 13(5-6):585--603,
  2007.

\bibitem{aub}
Thierry Aubin.
\newblock Meilleures constantes dans le th{\'e}or{\`e}me d'inclusion de
  {S}obolev et un th{\'e}or{\`e}me de {F}redholm non lin{\'e}aire pour la
  transformation conforme de la courbure scalaire.
\newblock {\em J. Funct. Anal.}, 32(2):148--174, 1979.

\bibitem{bdm}
Daniele Bartolucci, Francesca De~Marchis, and Andrea Malchiodi.
\newblock Supercritical conformal metrics on surfaces with conical
  singularities.
\newblock {\em Int. Math. Res. Not. IMRN}, (24):5625--5643, 2011.

\bibitem{bm}
Daniele Bartolucci and Andrea Malchiodi.
\newblock An improved geometric inequality via vanishing moments, with
  applications to singular liouville equations.
\newblock {\em Comm. Math. Phys., to appear, d.o.i. 10.1007/s00220-013-1731-0.}

\bibitem{bjmr}
Luca Battaglia, Aleks Jevnikar, Andrea Malchiodi, and David Ruiz.
\newblock A general existence result for the toda system on compact surfaces,
  http://arxiv.org/abs/1306.5404.

\bibitem{breme}
Haim Brezis and Frank Merle.
\newblock Uniform estimates and blow-up behavior for solutions of {$-\Delta
  u=V(x)e^u$} in two dimensions.
\newblock {\em Comm. Partial Differential Equations}, 16(8-9):1223--1253, 1991.

\bibitem{cay}
Luis~A. Caffarelli and Yi~Song Yang.
\newblock Vortex condensation in the {C}hern-{S}imons {H}iggs model: an
  existence theorem.
\newblock {\em Comm. Math. Phys.}, 168(2):321--336, 1995.

\bibitem{carma}
Alessandro Carlotto and Andrea Malchiodi.
\newblock Weighted barycentric sets and singular {L}iouville equations on
  compact surfaces.
\newblock {\em J. Funct. Anal.}, 262(2):409--450, 2012.

\bibitem{ci}
Dongho Chae and Oleg~Yu. Imanuvilov.
\newblock The existence of non-topological multivortex solutions in the
  relativistic self-dual {C}hern-{S}imons theory.
\newblock {\em Comm. Math. Phys.}, 215(1):119--142, 2000.

\bibitem{cysing}
Sun-Yung~A. Chang and Paul~C. Yang.
\newblock Conformal deformation of metrics on {$S^2$}.
\newblock {\em J. Differential Geom.}, 27(2):259--296, 1988.

\bibitem{cy0}
Sun-Yung~Alice Chang and Paul~C. Yang.
\newblock Prescribing {G}aussian curvature on {$S^2$}.
\newblock {\em Acta Math.}, 159(3-4):215--259, 1987.

\bibitem{clin}
Chiun-Chuan Chen and Chang-Shou Lin.
\newblock Topological degree for a mean field equation on {R}iemann surfaces.
\newblock {\em Comm. Pure Appl. Math.}, 56(12):1667--1727, 2003.

\bibitem{cl3}
Chiun-Chuan Chen and Chang-Shou Lin.
\newblock Mean field equations of {L}iouville type with singular data: sharper
  estimates.
\newblock {\em Discrete Contin. Dyn. Syst.}, 28(3):1237--1272, 2010.

\bibitem{chenwx}
Wen~Xiong Chen.
\newblock A {T}r{\"u}dinger inequality on surfaces with conical singularities.
\newblock {\em Proc. Amer. Math. Soc.}, 108(3):821--832, 1990.

\bibitem{chen}
Xiuxiong Chen.
\newblock Remarks on the existence of branch bubbles on the blowup analysis of
  equation {$-\Delta u=e^{2u}$} in dimension two.
\newblock {\em Comm. Anal. Geom.}, 7(2):295--302, 1999.

\bibitem{csw}
M.~Chipot, I.~Shafrir, and G.~Wolansky.
\newblock On the solutions of {L}iouville systems.
\newblock {\em J. Differential Equations}, 140(1):59--105, 1997.

\bibitem{cia}
Andrea Cianchi.
\newblock Moser-{T}rudinger trace inequalities.
\newblock {\em Adv. Math.}, 217(5):2005--2044, 2008.

\bibitem{djlw}
Weiyue Ding, J{\"u}rgen Jost, Jiayu Li, and Guofang Wang.
\newblock Existence results for mean field equations.
\newblock {\em Ann. Inst. H. Poincar\'e Anal. Non Lin\'eaire}, 16(5):653--666,
  1999.

\bibitem{dja}
Zindine Djadli.
\newblock Existence result for the mean field problem on {R}iemann surfaces of
  all genuses.
\newblock {\em Commun. Contemp. Math.}, 10(2):205--220, 2008.

\bibitem{dunne1995self}
G.~Dunne.
\newblock {\em Self-dual Chern-Simons Theories}.
\newblock Lecture notes in physics. New series m: Monographs. Springer, 1995.

\bibitem{fm}
Luigi Fontana and Carlo Morpurgo.
\newblock Adams inequalities on measure spaces.
\newblock {\em Adv. Math.}, 226(6):5066--5119, 2011.

\bibitem{jlw}
J{\"u}rgen Jost, Changshou Lin, and Guofang Wang.
\newblock Analytic aspects of the {T}oda system. {II}. {B}ubbling behavior and
  existence of solutions.
\newblock {\em Comm. Pure Appl. Math.}, 59(4):526--558, 2006.

\bibitem{jw}
J{\"u}rgen Jost and Guofang Wang.
\newblock Analytic aspects of the {T}oda system. {I}. {A} {M}oser-{T}rudinger
  inequality.
\newblock {\em Comm. Pure Appl. Math.}, 54(11):1289--1319, 2001.

\bibitem{kaolee}
Hsien-Chung Kao and Kimyeong Lee.
\newblock Self-dual {${\rm SU}(3)$} {C}hern-{S}imons {H}iggs systems.
\newblock {\em Phys. Rev. D (3)}, 50(10):6626--6632, 1994.

\bibitem{lee}
Kimyeong Lee.
\newblock Self-dual nonabelian {C}hern-{S}imons solitons.
\newblock {\em Phys. Rev. Lett.}, 66(5):553--555, 1991.

\bibitem{linweizhang}
Changshou Lin, Juncheng Wei, and Lei Zhang.
\newblock Classification of blowup limits for su(3) singular toda systems.
\newblock {\em preprint}, 2013.

\bibitem{cheikh}
Andrea Malchiodi and Cheikh~Birahim Ndiaye.
\newblock Some existence results for the {T}oda system on closed surfaces.
\newblock {\em Atti Accad. Naz. Lincei Cl. Sci. Fis. Mat. Natur. Rend. Lincei
  (9) Mat. Appl.}, 18(4):391--412, 2007.

\bibitem{malrui}
Andrea Malchiodi and David Ruiz.
\newblock New improved {M}oser-{T}rudinger inequalities and singular
  {L}iouville equations on compact surfaces.
\newblock {\em Geom. Funct. Anal.}, 21(5):1196--1217, 2011.

\bibitem{malruicpam}
Andrea Malchiodi and David Ruiz.
\newblock A variational analysis of the {T}oda system on compact surfaces.
\newblock {\em Comm. Pure Appl. Math.}, 66(3):332--371, 2013.

\bibitem{mos}
J.~Moser.
\newblock A sharp form of an inequality by {N}. {T}rudinger.
\newblock {\em Indiana Univ. Math. J.}, 20:1077--1092, 1970/71.

\bibitem{ntcv99}
Margherita Nolasco and Gabriella Tarantello.
\newblock Double vortex condensates in the {C}hern-{S}imons-{H}iggs theory.
\newblock {\em Calc. Var. Partial Differential Equations}, 9(1):31--94, 1999.

\bibitem{noltar2000}
Margherita Nolasco and Gabriella Tarantello.
\newblock Vortex condensates for the {${\rm SU}(3)$} {C}hern-{S}imons theory.
\newblock {\em Comm. Math. Phys.}, 213(3):599--639, 2000.

\bibitem{sy}
Joel Spruck and Yi~Song Yang.
\newblock Topological solutions in the self-dual {C}hern-{S}imons theory:
  existence and approximation.
\newblock {\em Ann. Inst. H. Poincar{\'e} Anal. Non Lin{\'e}aire},
  12(1):75--97, 1995.

\bibitem{tar96}
Gabriella Tarantello.
\newblock Multiple condensate solutions for the {C}hern-{S}imons-{H}iggs
  theory.
\newblock {\em J. Math. Phys.}, 37(8):3769--3796, 1996.

\bibitem{tarl}
Gabriella Tarantello.
\newblock {\em Selfdual gauge field vortices}.
\newblock Progress in Nonlinear Differential Equations and their Applications,
  72. Birkh{\"a}user Boston Inc., Boston, MA, 2008.
\newblock An analytical approach.

\bibitem{tro}
Marc Troyanov.
\newblock Prescribing curvature on compact surfaces with conical singularities.
\newblock {\em Trans. Amer. Math. Soc.}, 324(2):793--821, 1991.

\bibitem{tru}
Neil~S. Trudinger.
\newblock On imbeddings into {O}rlicz spaces and some applications.
\newblock {\em J. Math. Mech.}, 17:473--483, 1967.

\bibitem{ty}
Anatoly Tur and Vladimir Yanovsky.
\newblock Point vortices with a rational necklace: new exact stationary
  solutions of the two-dimensional {E}uler equation.
\newblock {\em Phys. Fluids}, 16(8):2877--2885, 2004.

\bibitem{wang99}
Guofang Wang.
\newblock Moser-{T}rudinger inequalities and {L}iouville systems.
\newblock {\em C. R. Acad. Sci. Paris S{\'e}r. I Math.}, 328(10):895--900,
  1999.

\bibitem{yyang}
Yisong Yang.
\newblock {\em Solitons in field theory and nonlinear analysis}.
\newblock Springer Monographs in Mathematics. Springer-Verlag, New York, 2001.

\end{thebibliography}
\bibliographystyle{plain}

\end{document}